\documentclass[12pt,leqno]{amsart}
\usepackage[dvips]{graphicx}
\usepackage{amsfonts}
\usepackage{amssymb}
\usepackage{amsmath}
\usepackage{amscd}
\usepackage{graphicx}
\def\DATE{\today}
\newtheorem{theorem}{Theorem}
\newtheorem{definition}[theorem]{Definition}

\newtheorem{lemma}[theorem]{Lemma}

\newtheorem{proposition}[theorem]{Proposition}

\newcommand\g{\mathfrak{g}}

\newcommand\K{\mathbb{K}}

\newcommand\lr{\left\{ \begin{array}{l}}
\def\ds{\displaystyle}

\title{Breadth and characteristic sequence of nilpotent Lie algebras }
\author{Elisabeth Remm}
\date{Chevat 16, 5775}
\address{Universit\'{e} de Haute Alsace, LMIA, 6 rue des Fr\`{e}res Lumi\`{e}re, 68093 Mulhouse}
\email{elisabeth.remm@uha.fr}

\begin{document}

\maketitle

\begin{abstract}
The notion of breadth of a nilpotent Lie algebra was introduced and used  to approach problems of classification up to  isomorphism in \cite{KMS}. In the present paper, we study this invariant in terms of characteristic sequence, another invariant introduced by Goze and Ancochea in \cite{AG1}. This permits to complete the determination of Lie algebras of breadth $2$ studied in \cite{KMS} and to begin the work for Lie algebras with breadth greater than $2$.

\end{abstract} 

\noindent{\bf Introduction}

The breadth of a Lie algebra was introduced in \cite{KMS}. As it is well explained in this paper, initially, the notion of breadth was defined to investigate the $p$-groups: the breadth of a $p$-group is the size of the largest conjugacy class in this $p$-group. 
Recently  characterizations of finite p-groups of breadth 1 and 2 have been done.
 In \cite{KMS}, the authors translate this notion to finite dimensional nilpotent Lie algebras: the breadth of a nilpotent Lie algebra is the maximum of the ranks of the adjoint operators and in their paper, they describe nilpotent Lie algebras of breadth $1$ and, for dimensions $5$ and $6$, of breadth $2$. 

We know that the classification, up to isomorphism, of finite dimensional nilpotent Lie algebras is a very difficult problem. It is solved up to the dimension $7$ in the complex and real cases and  for some special families of nilpotent Lie algebras in bigger dimensions. In  recent papers \cite{GRkstep,RFili}, since the classification problem in large dimensions  seems unrealistic, we choose to describe classes of nilpotent Lie algebras, in any dimension, associated with the characteristic sequence. For these families, we have reduced the notion of Chevalley-Eilenberg cohomology and given, sometimes, topological properties of these families. In this paper, we consider these results to complete in any dimension the characterization of nilpotent Lie algebras of breadth equal to $2$, completing, in some way, the work \cite{KMS}. We develop also the case where the breadth is $3$.

\section{breadth and characteristic sequence}
\subsection{Breadth of a nilpotent Lie algebra}
Let $\g$ be a finite dimensional Lie algebra over field $\K$ of characteristic $0$.
\begin{definition} \cite{KMS}
The breadth of $\g$ is the invariant, up to isomorphism,
$$b(\g)=\max \{{\rm rank}({\rm ad} X), \ X \in \g\}.$$
More generally, if $I$ is a not trivial ideal of $\g$, then the breadth of $\g$  on $I$ is
$$b_I(\g)=\max \{ {\rm rank}({\rm ad} X|_I), \ X \in \g\}.$$
\end{definition} 
We have $b_{\g}(\g)=b(\g)$ and clearly $b_I(\g) \leq b(\g)$. For example, $\g$ is abelian if and only if $b(\g) =0$ and $\g$ is filiform (see \cite{GK}) if and only if $b(\g)=\dim \g - 2.$ Recall that a finite dimensional nilpotent Lie algebra is called filiform, notion introduced by Michele Vergne, if $\g$ satisfies
$$\dim \mathcal{C}_i(\g)=i$$
for $i=0,\cdots,n-2$ where $n=\dim \g$ and  $\mathcal{C}_i(\g)$ are the terms of the ascending sequence of  $\g$ that is 
$$\left\{ 
\begin{array}{l}
 \mathcal{C}_0(\g)=\{ 0 \}, \\
  \mathcal{C}_i(\g)= \{ X \in \g, [X,\g]\in  \mathcal{C}_{i-1}(\g)\} .
\end{array}
\right.$$
 The following proposition is proved in \cite{KMS}
\begin{proposition}\label{lowerbound}
For any finite dimensional Lie algebra $\g$, we have $$b(\g) \leq {\rm dim}(\g/Z(\g)) -1.$$
\end {proposition}
In the case of a filiform Lie algebra, the upper bound of this inequality is satisfied. In fact  if $\g$ is filiform, then $\dim Z(\g)=1$ and $b(\g)=n-2= {\rm dim}(\g/Z(\g)) -1.$  

Assume now that $\g$ is finite dimensional and nilpotent.  We have seen that $b(\g)=0$ if and only if $\g$  is abelian. Now if $b(\g)=1$, then, from \cite{KMS} we have $\dim  \mathcal{C}^1(\g)=1$ where $\mathcal{C}^1(\g)$ is the derived subalgebra of $\g,$ that is generated by the Lie brackets $[X,Y]$ for any $X,Y \in \g$. Then, if $\g$ is indecomposable, that is without direct factor, then $\g$ is an Heisenberg Lie algebra. The non trivial study begins with $b(\g)=2$. This is a large part of the work of \cite{KMS}. We will complete their result using and comparing with an another invariant of nilpotent Lie algebras, the characteristic sequence, presented in the following section.

\subsection{Characteristic sequence}
Let $\g$ be a finite dimensional nilpotent Lie algebra over  $\K$. Let $\mathcal{C}^1(\g)$ and $Z(\g)$ be respectively the derived subalgebra  and  center of $\g$. The characteristic sequence of $\g$, sometimes called the Goze's invariant and  introduced in \cite{AG1},  is an useful invariant in the classification theory of finite-dimensional nilpotent Lie algebras. It is define as follows: for any $X \in \g,$ let $c(X)$ be the ordered sequence, for the lexicographic order, of the dimensions of the Jordan blocks of the nilpotent operator $ad X$. The characteristic sequence of $\g$ is the invariant, up to isomorphism,
$$c(\g)=\max \{c(X), \ X \in \g-\mathcal{C}^1(g)\}.$$
Then $c(\g)$ is a sequence of type $(c_1,c_2,\cdots,c_k,1)$ with $c_1 \geq c_2 \geq \cdots \geq c_k \geq 1$ and $c_1+c_2+\cdots+c_k+1= n =\dim \g$. For example, \begin{enumerate}
  \item $c(\g)=(1,1,\cdots,1)$ if and only if $\g$ is abelian,
  \item  $c(\g)=(2,1,\cdots,1)$ if and only if $\g$ is a direct product of an Heisenberg Lie algebra by an abelian ideal,
 \item  If $\g$ is $2$-step nilpotent then there exists $p$ and $q$ such that $c(\g)=(2,2,\cdots,2,1,\cdots,1)$ with $2p+q=n$, that is $p$ is the occurrence of $2$ in the characteristic sequence and $q$ the occurrence of $1$.
 \item $\g$ is filiform if and only if $c(\g)=(n-1,1)$.
 \end{enumerate}
 A vector $X \in \g-\mathcal{C}^1(\g)$ is called a characteristic vector if $c(X)=c(\g)$. Of course, such a vector is not unique.

\subsection{Breadth and characteristic sequence}

Let $\g$ be a finite dimensional nilpotent Lie algebra, $b(\g)$ its breadth and $c(\g)$ its characteristic sequence. 
\begin{lemma}
Let $X$ be a characteristic vector with $c(\g)=c(X)=(c_1,\cdots,c_k,1)$. Then
\begin{equation}
\label{ rank}
{\rm rank}({\rm ad} X)=\ds \sum_{i=1}^k c_i-k
\end{equation}
\end{lemma}
In fact, if we consider the Jordan block $J_i$ of dimension $c_i$ of the operator ${\rm ad} X$, since $0$ is the only eigenvalue, we  have ${\rm rank}(J_i)=c_i-1.$

We deduce

\begin{theorem} \label{b}
Let $\g$ be a finite dimensional Lie algebra whose characteristic sequence is $c(\g)=(c_1,\cdots,c_k,1)$. Then, its breadth satisfies
\begin{equation}
\label{ breadth}
b(\g)=\ds \sum_{i=1}^k c_i-k.
\end{equation}
\end{theorem}

\medskip

\noindent{\bf Consequences.}

1. If $\g$ is $n$-dimensional $2$-step nilpotent, then $c(\g)=(2,\cdots,2,1,\cdots,1)$. If $2$ appears $p$ times, then
$$b(\g)=p.$$

2.  If $\g$ is $n$-dimensional $3$-step nilpotent, then $c(\g)=(3,\cdots,3,2,\cdots,2,1,\cdots,1)$. If $3$ appears $p_3$ times, $2$ appears $p_2$ times in this sequence, then
$$b(\g)=2p_3+p_2.$$

3. More generally, if $\g$ is $n$-dimensional $k$-step nilpotent, then  its characteristic sequence is  $c(\g)=(k,\cdots,k,\cdots ,3,\cdots,3,2,\cdots,2,1,\cdots,1)$. Assume that $k$ appears $p_k$-times, $k-1$ $p_{k-1}$-times, up to $2$ which appears $p_{2}$-times and $1$ which appears $p_1$-times. Then
$$b(\g)=(k-1)p_k+(k-2)p_{k-1}+\cdots+2p_{3}+p_{2}=\sum_{i=1}^{k}(i-1)p_i.$$

\begin{proposition}\label{b=2}
A $n$-dimensional nilpotent Lie algebra satisfies $b(\g)=2$ if and only if one of these conditions is satisfied:\begin{enumerate}
  \item $\g$ is $2$-step nilpotent with $c(\g)=(2,2,1,\cdots,1)$.
  \item $\g$ is $3$-step nilpotent with $c(\g)=(3,1,\cdots,1)$.
\end{enumerate}
\end{proposition}
In fact, if $\g$ is $k$-step nilpotent with $k \geq 4$, then $b(\g)=(k-1)p_k+(k-2)p_{k-1}+\cdots+2p_{3}+p_{2}$ with $p_k >0$. Then $b(\g) \geq 3.$ Now if $\g$ is $3$-step nilpotent, then $b(\g)=2p_3+p_2$ with $p_3 >0$. Then $b(\g)=2$ implies $p_3=1$ and $p_2=0$. If $\g$ is $2$-step nilpotent, then $b(\g)=p_2$ this gives $p_2=2$.

\section{Nilpotent Lie algebras with breadth equal to $2$}

 A lot of results concerning nilpotent Lie algebras with breadth equal to $2$ are presented in \cite{KMS}, principally for the small dimensions. In this section, we propose to complete this study using Proposition \ref{b=2}. But, before  begining this study, we recall some properties of an important tool used in \cite{GRkstep} and based on  deformation theory of Lie algebras.
Let $\g_0= (\K^n,\mu_0)$ be a $n$-dimensional  Lie algebra over $\K$ where $\mu_0$ denotes the Lie bracket of $\g_0$ and $\K^n$ the underlying linear space to $\g_0$. A formal deformation $\g$ of $\g_0$ is given by a pair
$\g=(\K^n,\mu)$ where $\mu$ is a formal series
$$\mu=\mu_0 + \sum t^i \varphi_i$$
with  skew-symmetric bilinear maps $\varphi_i$ on $\K^n$ and $\mu$ satisfying the Jacobi identity.  To simplify notations, we denote for any bilinear map $\varphi$ on $\K^n$ with values on $\K^n$,  $\varphi \bullet \varphi$ the trilinear map
$$ \varphi \bullet \varphi(X,Y,Z)= \varphi (  \varphi(X,Y),Z)+\varphi (  \varphi(Y,Z),X)+\varphi (  \varphi(Z,X),Y)$$
and by $\varphi\circ_1 \varphi$ the ${\rm comp}_1$-operation (the composition on the first element), that is
$$\varphi\circ_1 \varphi (X,Y,Z)=\varphi (  \varphi(X,Y),Z).$$
Then the Jacobi identity for $\varphi$ is equivalent to
$\varphi \bullet \varphi =0.$ If $\mu$ is a formal deformation of $\mu_0$, then $\mu \bullet \mu =0$ implies
$$
\left\{
\begin{array}{l}
   \delta_{C, \mu_0}\varphi_1=0 ,    \\
     \varphi_1 \bullet \varphi_1=   \delta_{C, \mu_0}\varphi_2,
\end{array}
\right.
$$
where $ \delta_{C, \mu_0}$ is the coboundary operator associated to the Chevalley-Eilenberg complex of $\g_0$:
$$ \delta_{C, \mu_0}(\varphi)=\mu_0 \bullet \varphi + \varphi \bullet \mu_0.$$ We shall denote by $H^*_{C}(\g_0,\g_0)$ or sometimes
$H^*_{C}(\mu_0,\mu_0)$ the corresponding cohomology spaces. A special class of formal deformation is the class of linear deformations. We call {\it linear deformation} of $\g_0= (\K^n,\mu_0)$ a formal deformation $\mu=\mu_0 + \sum t^i \varphi_i$ with $\varphi_k=0$ for any $k \geq 2$, that is 
$$\mu=\mu_0+t\varphi_1.$$
In this case, the Jacobi identity for $\mu$ is equivalent to
\begin{equation}
\label{line}
\left\{
\begin{array}{l}
   \delta_{C, \mu_0}\varphi_1=0 ,    \\
     \varphi_1 \bullet \varphi_1= 0.
     \end{array}
\right.
\end{equation}

\subsection{Lie algebras with characteristic sequence $(3,1,\cdots,1)$}
Let $\g=(\K^n,\mu)$ be a $n$-dimensional nilpotent Lie algebra with characteristic sequence $c(\g)=(3,1,\cdots,1)$. Then from \cite{GRkstep} it is isomorphic to a linear deformation of 
the Lie algebra $\g_{1,0,n-1}=(\K^n,\mu_0)$ whose Lie brackets are
$$[X_1,X_2]=X_3, \ [X_1,X_3]=X_4, \ [X_1,X_i]=0 \ {\rm for} \  i =4, \cdots , n,  \ [X_i,X_j]=0, \ 2\leq i < j \leq n.$$
So $\mu=\mu_0+t\varphi$ where $\varphi$ is a bilinear map satisfying the conditions  (\ref{line}) and conditions which imply that $\mu$ is a Lie bracket of a $3$-step nilpotent Lie algebra. All these conditions are described in \cite{GRkstep} and imply that $\g$ is in  one of the following cases:
\begin{enumerate}
\item $\g$ is  isomorphic to $\g_{1,0,n-1}$,
  \item $\g$ is a direct product of an abelian ideal and of the $5$-dimensional Lie algebra and its Lie bracket is isomorphic to $\mu_0+\varphi$ with $\varphi(X_2,X_3)=aX_5, $
 \item $\g$ is a direct product of an abelian ideal and of the $5$-dimensional Lie algebra and its Lie bracket is isomorphic to $\mu_0+\varphi$ with $\varphi(X_2,X_5)=bX_4$ 
  \item If $\g$ is indecomposable and of dimension greater than $5$, then 
 $\g$ is isomorphic to the Lie algebra
$$
\left\{
\begin{array}{l}
[X_1,X_{2}\rbrack=X_{3}, \
\lbrack X_1,X_{3}\rbrack=X_{4},  \\
\lbrack X_2,X_{k}\rbrack=a_{2,k}X_4, \ k \geq 5,\\
\lbrack X_l,X_{k}\rbrack=a_{l,k}X_{4}, \ 5 \leq l < k\leq n.
\end{array}
\right.
$$
\end{enumerate} 
This last family of Lie algebras, parametrized by the constants structure $a_{2,k}$ and $a_{l,k}$ can be reduced to obtain a classification up to isomorphism. 
\begin{proposition}  \cite{CGJ}
Any $n$-dimensional nilpotent Lie algebra, $n\geq 4$ with breadth equal to $2$ and characteristic sequence  $(3,1,\cdots,1)$ is isomorphic to one of the following Lie algebras whose Lie brackets are
\begin{enumerate}
  \item $ [X_1,X_2]=X_3, \ [X_1,X_3]=X_4$ {\rm that is} \  $\g_{1,0,n-1}$.
    \item  $ [X_1,X_2]=X_3, \ [X_1,X_3]=X_4, [X_2,X_3]=X_5,$  {\rm here} $n \geq 5$.
  \item  $ [X_1,X_2]=X_3, \ [X_1,X_3]=X_4,[X_{2i+1},X_{2i+2}]=X_4, i=1,\cdots,p-1$  {\rm  with}  $n=2p$,
  \item  $ [X_1,X_2]=X_3, \ [X_1,X_3]=X_4,[X_2,X_n]=X_4,[X_{2i+1},X_{2i+2}]=X_4, i=1,\cdots,p-1$  {\rm  with} $n=2p+1$,
\end{enumerate} 
 \end{proposition}

 \subsection{Lie algebras whose characteristic sequence is $(2,2,1,\cdots,1)$}
 
 In \cite{GRsymplectique}, we have studied nilpotent Lie algebras with characteristic sequence $(2,2,1,1,1,1)$ in order to construct symplectic form on these Lie algebras. We will use notations and results which are in this paper.
Let $\mathcal{F}^{n,2}$ be the set of $n$-dimensional $2$-step nilpotent Lie algebras and 
 $\mathcal{F}^{n,2}_{2,n-2}$  the subset  corresponding to Lie algebras of characteristic sequence $(2,2,1,\cdots,1)$ (the subscripts correspond to the numbers of  $2$ and $1$ in the characteristic sequence). Let $\g_{2,n-2}$  be the $n$-dimensional Lie algebra given by the brackets
$$
[X_1,X_{2i}]=X_{2i+1}, \quad i=1,2.
$$
It is obvious that $\g_{2,n-2}
 \in \mathcal{F}^{n,2}_{2,n-2}$. From \cite{GRsymplectique}, any  $n$-dimensional  $2$-step nilpotent Lie algebra with characteristic sequence $(2,2,1,\cdots,1)$  is isomorphic to a linear deformation of $\g_{2,n-2}$. This means that its Lie bracket is isomorphic to $\mu=\mu_0 + t \varphi$,  where $\mu_0$ is the Lie bracket of $\g_{2,n-2}$, and $\varphi$ is a skew-bilinear map such that
$$
\left\{
\begin{aligned}
&    \varphi \in Z^2_{CH}(\g_{2,n-2},\g_{2,n-2}), \\
&     \varphi \circ_1 \varphi=0.
\end{aligned}
\right.
$$
Recall that $Z^2_{CH}(\g_{2,n-2},\g_{2,n-2})$ is constituted of bilinear maps $\varphi$ on $\K^n$ with values in this space such that
$$\varphi \circ_1\mu_0+\mu_0\circ_1 \varphi=0.$$
In fact, in \cite{GRkstep}, we have define a sub-complex of the Chevalley-Eilenberg complex whose  main property is that for any deformation $\mu=\mu_0 + \sum t^i\varphi_i$ of a $2$-step nilpotent Lie algebra $\g_0=(\mu_0,\K^n)$ which is also $2$-step nilpotent, then $\varphi_1 \in Z^2_{CH}(\mu_0,\mu_0)$. 

We deduce that
the family $\mathcal{F}^{n,2}_{2,n-2}$  of  $n$-dimensional $2$-step nilpotentLie algebras with characteristic sequence $(2,2,1, \cdots,1)$
 is the union of two algebraic components, the first one, $\mathcal{C}_1(\mathcal{F}^{n,2}_{2,n-2})$, corresponding to the cocycles
\begin{equation}
\label{C3}
\varphi(X_{i},X_{j})= \sum_{k=1}^{2}a_{i,j}^{2k+1}X_{2k+1}, \quad  2\leq i<j \leq n, \quad i,j \neq 3,5
\end{equation}
the second one, $\mathcal{C}_2(\mathcal{F}^{n,2}_{2,n-2})$,  to the cocyles
\begin{equation}
\label{C_2}
\left\{
\begin{array}{l}
\varphi(X_{i},X_{j})= \sum_{k=1}^{2}a_{i,j}^{2k+1}X_{2k+1},  \quad (i,j) \neq (2,4), \ 2\leq i < j \leq n-1,\ \  i,j \notin \{3,5\}\\
\varphi(X_{2},X_{4})=X_n
\end{array}
\right.
\end{equation}
where the undefined products $\varphi(X,Y)$ are nul or obtained by skew-symmetry.
Each of these components is a  regular algebraic variety. These  components can be characterized by  the number of generators.
\begin{proposition}
Any $n$-dimensional nilpotent Lie algebra with breadth equal to $2$ and characteristic sequence $(2,2,1,\cdots,1)$ is isomorphic to a Lie algebra belonging to one of the components $\mathcal{C}_1(\mathcal{F}^{n,2}_{2,n-2})$ or $\mathcal{C}_2(\mathcal{F}^{n,2}_{2,n-2})$.
\end{proposition}

\section{Nilpotent Lie algebras with breadth equal to $3$}

From Theorem \ref{b}, the characteristic sequence of a $n$-dimensional Lie algebra of breadth $b(\g)=3$ is equal to one of the following:
\begin{enumerate}
  \item $c(\g)=(4,1,\cdots,1),$ ($n \geq 5$),
  \item $c(\g)=(3,2,1,\cdots,1)$, ($n \geq 6$),
  \item $c(\g)=(2,2,2,1,\cdots,1)$, ($n \geq 7$).
\end{enumerate}

\subsection{Nilpotent Lie algebras of characteristic sequence $(4,1,\cdots,1)$}

Let us consider the $n$-dimensional nilpotent Lie algebra $\g_{1,0,0,n-4}=(\K^n,\mu_0)$ with characteristic sequence $(4,1,\cdots,1)$ whose Lie brackets in the basis $\{X_1,\cdots,X_n\}$ are
$$ \mu_0(X_1,X_i)=X_{i+1}, \ i=2,3,4$$
all the others brackets being nul.
 Any $n$-dimensional Lie algebra $\g$ whose characteristic sequence is also $(4,1,\cdots,1)$ is isomorphic to a Lie algebra whose Lie bracket $\mu$ is a linear deformation of $\mu_0$, that is $\mu=\mu_0+t\varphi$
where $\varphi$ is a $2$-cocycle of the Chevalley-Eilenberg complex $H_C^*(\mu_0,\mu_0)$ satisfying 
\begin{equation}
\label{4step}
\left\{
\begin{array}{l}
\delta_{C,\mu_0}^2(\varphi)=0,\\
\varphi \bullet \varphi =0,\\
\delta_{R,\mu_0}^2(\varphi)=0,\\
\mu_0^2 \circ_1\varphi ^2+  \mu_0  \circ_1 \varphi \circ_1 \mu_0 \circ_1 \varphi + \mu_0\circ_1 \varphi ^2\circ_1 \mu_0 +\varphi\circ_1 \mu_0^2   \circ_1 \varphi +\\ +\varphi\circ_1 \mu_0  \circ_1 \varphi \circ_1 \mu_0+\varphi2\circ_1 \mu_0^2=0, \\
\mu_0 \circ_1\varphi ^3+\varphi\circ_1 \mu_0 \circ_1 \varphi^2  +\varphi^2\circ_1 \mu_0  \circ_1 \varphi+\varphi^3\circ_1\mu_0=0, \\
\varphi  ^4=0,
\end{array}
\right.
\end{equation}
where $\delta_{R,\mu_0}^2(\varphi)$ is the $5$-linear map
$$\delta_{R,\mu_0}^2(\varphi)=\mu_0 ^3\circ_1\varphi + \mu_0 ^2\circ_1 \varphi \circ_1 \mu_0 +\mu_0 \circ_1 \varphi \circ_1 \mu_0 ^2+\varphi\circ_1 \mu_0 ^ 3,$$
 and $\varphi^k=\varphi\circ_1\cdots\circ_1 \varphi.$  Moreover, $\varphi$ satisfies $\varphi(X_1,Y)=0$ for any $Y$ and $X_1$ is not in the derived subalgebra of $\g$. This implies that $(\mu_0 \circ_1 \varphi )(X_i,X_j,X_k)=0 $ if $1 <i <j <k $. Then (\ref{4step}) gives
 \begin{equation}
\left\{
\begin{array}{l}
\varphi(X_5,X_i)=0, \ \forall i,\\
\varphi(X_2,X_3)=aX_5+\sum_{k=6}^nd_kX_k,\\
\varphi(X_2,X_k)=a_kX_4+b_kX_5 \ \ 6 \leq k \leq n, \\
\varphi(X_3,X_k)=a_kX_5  \ \ 6 \leq k \leq n, \\
\varphi (X_i,X_j)= c_{i,j}X_5 \ \  6 \leq i < j \leq n.
\end{array}
\right.
\end{equation}
\begin{proposition}
Every $n$-dimensional nilpotent Lie algebra ($n\geq 5$) of breadth $3$ and characteristic sequence $(4,1,\cdots,1)$ is isomorphic to a Lie algebra belonging to one of the families
$$
\left\{
\begin{array}{l}
[X_1,X_i]=X_{i+1}, \ \ i=2,3,4,\\
\lbrack X_5,X_i\rbrack =0, \ \forall i,\\
\lbrack X_2,X_3\rbrack =aX_5,\\
\lbrack X_2,X_k\rbrack =a_kX_4+b_kX_5, \ \ 6 \leq k \leq n, \\
\lbrack X_3,X_k\rbrack =a_kX_5,  \ \ 6 \leq k \leq n, \\
\lbrack X_i,X_j\rbrack = c_{i,j}X_5, \ \  6 \leq i < j \leq n.
\end{array}
\right.
$$
$$
\left\{
\begin{array}{l}
[X_1,X_i]=X_{i+1}, \ \ i=2,3,4,\\
\lbrack X_5,X_i\rbrack =0, \ \forall i\,\\
\lbrack X_2,X_3\rbrack =aX_5+\sum_{k=6}^nd_kX_k,\\
\lbrack X_2,X_k\rbrack =b_kX_5, \ \ 6 \leq k \leq n, \\
\lbrack X_i,X_j\rbrack = c_{i,j}X_5, \ \  6 \leq i < j \leq n, \\
{\rm with} \ \displaystyle \sum_{i=6}^{j-1} c_{i,j}d_i-\sum_{i=j+1}^{n} c_{j,i}d_i =0 \ {\rm for} \  6 \leq  j \leq n. 
\end{array}
\right.
$$
\end{proposition}

The first family is a plane of dimension $\frac{(n-4)(n-3)}{2}.$  The second family defines an algebraic variety of codimension $n-5$ is a space of dimension $\frac{(n-4)(n-3)}{2}.$

Proposition 8 gives all the nilpotent Lie algebras of breadth $3$ and characteristic sequence $(4,1,\cdots,1).$ The nilpotent Lie algebras  of breadth $3$ and characteristic sequences $(3,2,1,\cdots,1)$ and $(2,2,2,1,\cdots,1)$ are determinated in \cite{GRkstep}.

\medskip

\noindent{\bf Remark.} We have just seen that the set $\mathcal{F}^{n,4}_{(1,0,0,n-5)}$ of $n$-dimensional Lie algebras, $n \geq 5$, whose characteristic sequence is $(4,1,\cdots,1)$ is an algebraic variety which is the union of two irreducible components, the first one is the orbit by the linear group of a $\frac{(n-4)(n-3)}{2}$ plane then it is reduced (the affine shema which corresponds to this component is reduced). Then, if there exists a rigid Lie algebra $\g_1$ in $\mathcal{F}^{n,4}_{(1,0,0,n-5)}$, that is any deformation of $\g_1$ in $\mathcal{F}^{n,4}_{(1,0,0,n-5)}$ is isomorphic to $\g_1$, then $\dim H^2_{CR}(\g_1,\g_1)=0$, where $H^2_{CR}(\g_1,\g_1)=\frac{Z^2_{CR}(\g_1,\g_1)}{B^2_{CR}(\g_1,\g_1)}$ and the two cocycles of $Z^2_{CR}(\g_1,\g_1)$ are bilinear maps satisfying (\ref{4step}) and $B^2_{CR}(\g_1,\g_1)$ the subspace $\{ \delta f, f \in gl(n,\K)\} \cap Z^2_{CR}(\g_1,\g_1)$. Now, if we consider the Lie algebra 
$$
\left\{
\begin{array}{l}
[X_1,X_i]=X_{i+1}, \ \ i=2,3,4,\\
\lbrack X_2,X_3\rbrack =X_5,\\
\lbrack X_2,X_{n-1}\rbrack =X_5,\\
\lbrack X_2,X_{n}\rbrack =X_4, \\
\lbrack X_3,X_n\rbrack =X_5, \\
\lbrack X_6,X_7\rbrack = \cdots=\lbrack X_{r},X_n\rbrack=X_5,  \  r=2p \ {\rm if} \  n=2p+1 \ {\rm or} \  r=2p-1 \   {\rm if} \  n=2p.
\end{array}
\right.
$$
we have $\dim Z^2_{CR}(\g_1,\g_1)=\frac{(n-4)(n-3)}{2}= \dim B^2_{CR}(\g_1,\g_1).$ We deduce that $\g_1$ is rigid in $\mathcal{F}^{n,4}_{(1,0,0,n-5)}$ and the first component is the Zariski closure of the orbit, for the natural action of the linear group, of $\g_1$. A similar computation for the second component shows that the Lie algebra $\g_2$ given by
$$
\left\{
\begin{array}{l}
[X_1,X_i]=X_{i+1}, \ \ i=2,3,4,\\
\lbrack X_2,X_3\rbrack =X_n,\\
\lbrack X_2,X_n\rbrack =X_5, \\
\lbrack X_{2k},X_{2k+1} \rbrack = X_5 \ \  3 \leq k  \leq \frac{n-1}{2}, \\
\end{array}
\right.
$$
belongs to the second component and it is rigid in this component. So any $n$-dimensional Lie algebras of characteristic sequence $(4,1,\cdots,1)$ is a contraction (or degeneration) of one of these two Lie algebras.

\subsection{Lie algebras of characteristic sequence  $(2,2,2,1,\cdots,1)$}

Let $n \geq 8$ and let  $\g_{3,n-6}$ be the $n$-dimensional Lie algebra given by the brackets
$$[X_1,X_{2i}]=X_{2i+1}, \ 1\leq i \leq 3,$$
and other non defined brackets are equal to zero.
Any $2$-step nilpotent $n$-dimensional Lie algebra with characteristic sequence $(2,2,2,1,\cdots,1)$  is isomorphic to a linear deformation of $\g_{3,n-6}$. Its Lie bracket is isomorphic to $\mu=\mu_0 + t \varphi$,  where $\mu_0$ is the Lie bracket of $\g_{3,n-6}$, and $\varphi$ is a skew-bilinear form such that 
$$
\label{def}
\left\{
\begin{array}{l} 
\medskip
    \varphi \in Z^2_{CH}(\g_{3,n-6},\g_{3,n-6}), \\
     \varphi\circ_1\varphi=0.
\end{array}
\right.
$$
From the construction of the linear deformation, we can assume that $\varphi(X_1,X)=0$ for any $X$. Now the first relation is equivalent to
$$
\left\{
\begin{array}{l}
\medskip
\varphi(X_{2i},X_{2j})=\displaystyle \sum_{k=1}^{3}a_{2i,2j}^{2k+1}X_{2k+1}+\sum_{k=8}^n a_{2i,2j}^{k}X_{k}, \ 1\leq i<j \leq 3, \\
\medskip
\varphi(X_{s},X_{l})=  \displaystyle \sum_{k=1}^{3}b_{s,l}^{2k+1}X_{2k+1}, \ s,l \geq 8.
\end{array}
\right.
$$
In this case, the identity $ \varphi\circ_1\varphi=0$ is satisfied and with such cocycles we describe, up to isomorphism, all the deformations of $\g_{3,n-6}$ which are $2$-step nilpotent.
But, we have to restrict these deformations only  to Lie algebras with   characteristic sequence equal to $(2,2,2,1,\cdots,1)$. This implies that the breadth is equal to $3$. This reduces the expression of the cocycles to
\begin{enumerate}
  \item  $n\geq 9$
   $$
\left\{
\begin{array}{l}
\medskip
\varphi_1(X_{2},X_{4})=\displaystyle \sum_{k=1}^{3}a_{2,4}^{2k+1}X_{2k+1}+X_8, \\
\medskip
\varphi_1(X_{2},X_{6})=\displaystyle \sum_{k=1}^{3}a_{2,6}^{2k+1}X_{2k+1}+X_9, \\
\medskip
\varphi_1(X_{4},X_{6})=\displaystyle \sum_{k=1}^{3}a_{4,6}^{2k+1}X_{2k+1},\\
\medskip
\varphi_1(X_{s},X_{l})=  0 \ s,l \geq 8.
\end{array}
\right.
$$
  \item  $n \geq 8$
  $$
\left\{
\begin{array}{l}
\medskip
\varphi_2(X_{2},X_{4})=\displaystyle \sum_{k=1}^{3}a_{2,4}^{2k+1}X_{2k+1}+X_8, \\
\medskip
\varphi_2(X_{2i},X_{2j})=\displaystyle \sum_{k=1}^{3}a_{2i,2j}^{2k+1}X_{2k+1},  \ 2\leq i<j \leq 3, \\
\medskip
\varphi_2(X_{s},X_{l})=  a_{s,l}X_5 ,\ s,l \geq 8.
\end{array}
\right.
$$
  \item  $n \geq 8$
  $$
\left\{
\begin{array}{l}
\medskip
\varphi_3(X_{2i},X_{2j})=\displaystyle \sum_{k=1}^{3}a_{2i,2j}^{2k+1}X_{2k+1}, \ 1\leq i<j \leq 3, \\
\medskip
\varphi_3(X_{s},X_{l})=  \sum_{k=1}^3a_{s,l}X_{2k+1}, \ s,l \geq 8.
\end{array}
\right.
$$
\end{enumerate}

\begin{proposition}
Any $n$-dimensional nilpotent Lie algebra whose breadth is $3$ and characteristic sequence $(2,2,2,1,\cdots,1)$ is isomorphic to a Lie algebra whose Lie bracket is $\mu_0+\varphi_i$ where $\varphi_i, \ i=1,2,3$ are defined above.
\end{proposition}

For the $7$-dimensional case, since the classification of $7$-dimensional nilpotent Lie algebras is known, from this list we find the following $7$-dimensional nilpotent Lie algebras of characteristic sequence $(2,2,2,1)$ (the notations are those of \cite{livreADL})

$
\frak{n}_7^{120} :
\begin{array}{l}
\lbrack X_1, X_i]=X_{i+1}, \ i=2,4,6, \ \
\lbrack X_2 , X_4 ]= X_7. \\
\end{array}
$

\smallskip

$
\frak{n}_7^{121} :
\begin{array}{l}
\lbrack X_1, X_i]=X_{i+1}. \ i=2,4,6 \ \\
\end{array}
$

\smallskip

$
\frak{n}_7^{122} :
\begin{array}{l}
\lbrack X_1, X_i]=X_{i+1}, \ i=2,4,6, \ \
\lbrack X_4 , X_6 ]= X_7. \\
\end{array}
$

\smallskip

$
\frak{n}_7^{123} :
\begin{array}{l}
\lbrack X_1, X_i]=X_{i+1}, \ i=2,4,6, \ \
\lbrack X_2, X_4 ]= X_5, \
\lbrack X_4 , X_6 ]= X_3. \\
\end{array}
$

\subsection{Lie algebras of characteristic sequence $(3,2,1,\cdots,1)$}
Consider now the case of characteristic sequence $(3,2,1,\cdots,1)$.Thus we assume that the dimension of $\g$ is greater than $6$. In dimension $6$ any Lie algebra with characteristic sequence $(3,2,1)$ is isomorphic to

\medskip

$\frak{n}^{11}_6$: $
\left[ X_1,X_i \right]=X_{i+1}, \ i=2,3,5, \
\left[ X_5,X_6 \right]=X_4 \
$

\smallskip

$\frak{n}^{12}_6$: $
\left[ X_1,X_i \right]=X_{i+1}, \ i=2,3,5, \
\left[ X_2,X_5 \right]=X_4 \
$

\smallskip

$\frak{n}^{13}_6$: $
\left[ X_1,X_i \right]=X_{i+1}, \ i=2,3,5, \
\left[ X_2,X_3 \right]=X_6 \
\left[ X_2,X_5 \right]=X_6 \
$

\smallskip

$\frak{n}^{14}_6$: $
\left[ X_1,X_i \right]=X_{i+1}, \ i=2,3,5, \
\left[ X_2,X_3 \right]=X_4-X_6 \
\left[ X_2,X_5 \right]=X_6 \
$

\smallskip

$\frak{n}^{15}_6$: $
\left[ X_1,X_i \right]=X_{i+1}, \ i=2,3,5, \
\left[ X_2,X_5 \right]=X_6 \
\left[ X_5,X_6 \right]=X_4 \
$

\smallskip

$\frak{n}^{16}_6$: $
\left[ X_1,X_i \right]=X_{i+1}, \ i=2,3,5 \
\left[ X_2,X_3 \right]=X_4 \
$

\smallskip

$\frak{n}^{17}_6$: $
\left[ X_1,X_i \right]=X_{i+1}, \ i=2,3,5 \
$
\medskip

 In dimension $7$, any indecomposable Lie algebra with characteristic sequence $(3,2,1,1)$ is isomorphic to

$
\frak{n}_7^{93} :
\begin{array}{l}
\lbrack X_1, X_i]=X_{i+1}, \ i=2,3,5, \ \
\lbrack X_2 , X_5 ]= X_7. \\
\end{array}
$

\smallskip

$
\frak{n}_7^{94} :
\begin{array}{l}
\lbrack X_1, X_i]=X_{i+1}, \ i=2,3,5, \ \
\lbrack X_2 , X_5 ]= X_4, \
\lbrack X_2 , X_3 ]= X_7 \\
\end{array}
$

\smallskip

$
\frak{n}_7^{95} :
\begin{array}{l}
\lbrack X_1, X_i]=X_{i+1}, \ i=2,3,5 \ \\
\lbrack X_2 , X_3 ]= X_7 \\
\end{array}
$

\smallskip

$
\frak{n}_7^{96} :
\begin{array}{l}
\lbrack X_1, X_i]=X_{i+1}, \ i=2,3,5, \ \
\lbrack X_3 , X_5 ]= -X_4, \
\lbrack X_2 , X_6 ]= X_4. \\
\end{array}
$

\smallskip

$
\frak{n}_7^{97} :
\left\{
\begin{array}{l}
\lbrack X_1, X_i]=X_{i+1}, \ i=2,3,5, \ \
\lbrack X_2 , X_6 ]= X_4, \
\lbrack X_3 , X_5 ]= -X_4, \\
\lbrack X_2 , X_5 ]= X_7. \\
\end{array}
\right.
$

\smallskip

$
\frak{n}_7^{98} :
\left\{
\begin{array}{l}
\lbrack X_1, X_i]=X_{i+1}, \ i=2,3,5, \ \
\lbrack X_2 , X_6 ]= X_4, \
\lbrack X_3 , X_5 ]= -X_4, \\
\lbrack X_2 , X_5 ]= X_7, \
\lbrack X_5 , X_6 ]= X_4. \\
\end{array}
\right.
$

\smallskip

$
\frak{n}_7^{99} :
\begin{array}{l}
\lbrack X_1, X_i]=X_{i+1}, \ i=2,3,5, \ \
\lbrack X_2 , X_7 ]= X_6, \
\lbrack X_2 , X_3 ]= X_4. \\
\end{array}
$

\smallskip

$
\frak{n}_7^{100} :
\begin{array}{l}
\lbrack X_1, X_i]=X_{i+1}, \ i=2,3,5, \ \
\lbrack X_2 , X_7 ]= X_4, \
\lbrack X_5 , X_7 ]= X_6. \\
\end{array}
$

\smallskip

$
\frak{n}_7^{101} :
\begin{array}{l}
\lbrack X_1, X_i]=X_{i+1}, \ i=2,3,5, \ \
\lbrack X_2 , X_7 ]= X_6, \
\lbrack X_5 , X_7 ]= X_4. \\
\end{array}
$

\smallskip

$
\frak{n}_7^{102} :
\begin{array}{l}
\lbrack X_1, X_i]=X_{i+1}, \ i=2,3,5, \ \
\lbrack X_5 , X_7 ]= X_4. \\
\end{array}
$

\smallskip

$
\frak{n}_7^{103} :
\begin{array}{l}
\lbrack X_1, X_i]=X_{i+1}, \ i=2,3,5, \ \
\lbrack X_2 , X_7 ]= X_4. \\
\end{array}
$

\smallskip

$
\frak{n}_7^{104} :
\begin{array}{l}
\lbrack X_1, X_i]=X_{i+1}, \ i=2,3,5, \ \
\lbrack X_5 , X_7 ]= X_4, \
\lbrack X_2 , X_3 ]= X_4. \\
\end{array}
$

\smallskip

$
\frak{n}_7^{105} :
\begin{array}{l}
\lbrack X_1, X_i]=X_{i+1}, \ i=2,3,5, \ \
\lbrack X_2 , X_7 ]= X_4, \
\lbrack X_2 , X_3 ]= X_4. \\
\end{array}
$

\smallskip

$
\frak{n}_7^{106} :
\begin{array}{l}
\lbrack X_1, X_i]=X_{i+1}, \ i=2,3,5, \ \
\lbrack X_2 , X_7 ]= X_4, \
\lbrack X_5 , X_6 ]= X_4. \\
\end{array}
$

\smallskip

$
\frak{n}_7^{107} :
\begin{array}{l}
\lbrack X_1, X_i]=X_{i+1}, \ i=2,3,5, \ \
\lbrack X_5 , X_7 ]= X_3, \
\lbrack X_6 , X_7 ]= X_4. \\
\end{array}
$

\smallskip

$
\frak{n}_7^{108} :
\begin{array}{l}
\lbrack X_1, X_i]=X_{i+1}, \ i=2,3,5, \ \
\lbrack X_2 , X_7 ]= X_6, \
\lbrack X_2 , X_3 ]= X_6. \\
\end{array}
$

\smallskip

$
\frak{n}_7^{109} :
\begin{array}{l}
\lbrack X_1, X_i]=X_{i+1}, \ i=2,3,5, \ \
\lbrack X_5 , X_7 ]= X_6, \
\lbrack X_2 , X_3 ]= X_6. \\
\end{array}
$

\smallskip

$
\frak{n}_7^{110} :
\begin{array}{l}
\lbrack X_1, X_i]=X_{i+1}, \ i=2,3,5, \ \
\lbrack X_2 , X_3 ]= X_6, \
\lbrack X_5 , X_7 ]= X_4.  \\
\end{array}
$

\smallskip

$
\frak{n}_7^{111} :
\begin{array}{l}
\lbrack X_1, X_i]=X_{i+1}, \ i=2,3,5, \ \
\lbrack X_2 , X_7 ]= X_6, \
\lbrack X_2 , X_5 ]= X_4. \\
\end{array}
$

\smallskip

$
\frak{n}_7^{112} :
\begin{array}{l}
\lbrack X_1, X_i]=X_{i+1}, \ i=2,3,5, \ \
\lbrack X_2 , X_3 ]= X_6, \
\lbrack X_2 , X_7 ]= X_4. \\
\end{array}
$

\smallskip

$
\frak{n}_7^{113} :
\begin{array}{l}
\lbrack X_1, X_i]=X_{i+1}, \ i=2,3,5, \ \
\lbrack X_5 , X_7 ]= X_6, \
\lbrack X_5 , X_6 ]= X_4. \\
\end{array}
$

\smallskip

$
\frak{n}_7^{114} :
\left\{
\begin{array}{l}
\lbrack X_1, X_i]=X_{i+1}, \ i=2,3,5, \ \
\lbrack X_2 , X_7 ]= X_4, \
\lbrack X_5 , X_6 ]= X_4, \\
\lbrack X_5 , X_7 ]= X_6. \\
\end{array}
\right.
$

\smallskip

$
\frak{n}_7^{115} :
\left\{
\begin{array}{l}
\lbrack X_1, X_i]=X_{i+1}, \ i=2,3,5, \ \
\lbrack X_2 , X_5 ]= X_4, \
\lbrack X_5 , X_7 ]= X_3, \\
\lbrack X_6 , X_7 ]= X_4. \\
\end{array}
\right.
 $

\smallskip

$
\frak{n}_7^{116} :
\left\{
\begin{array}{l}
\lbrack X_1, X_i]=X_{i+1}, \ i=2,3,5, \ \
\lbrack X_3 , X_5 ]= -X_4, \
\lbrack X_2 , X_6 ]= X_4,  \\
\lbrack X_5 , X_7 ]= -X_4. \\
\end{array}
\right.
$

\smallskip

$
\frak{n}_7^{117}(\alpha) :
\left\{
\begin{array}{l}
\lbrack X_1, X_i]=X_{i+1}, \ i=2,3,5, \ \
\lbrack X_2 , X_5 ]= X_7, \
\lbrack X_2 , X_7 ]= X_4, \\
\lbrack X_5 , X_6 ]= X_4, \
\lbrack X_5 , X_7 ]= \alpha X_4. \\
\end{array}
\right.
$

\smallskip

$
\frak{n}_7^{118} :
\left\{
\begin{array}{l}
\lbrack X_1, X_i]=X_{i+1}, \ i=2,3,5, \ \
\lbrack X_2 , X_5 ]= X_7, \
\lbrack X_2 , X_6 ]= X_4, \\
\lbrack X_3 , X_5 ]= -X_4, \
\lbrack X_5 , X_7 ]= -\frac{1}{4}X_4. \\
\end{array}
\right.
 $

\bigskip

Let us consider now the general case. We know that any $n$-dimensional nilpotent Lie algebra with characteristic sequence $(3,2,1,\cdots,1)$ is isomorphic to a linear deformation of $\g_{1,1,n-5}$ whose Lie bracket $\mu_0$ is given by
$$\mu_0(X_1,X_i)=X_{i+1}, \ \ i=2,3,5,$$
other non defined product are equal to zero. Since the classification up an isomorphism seems also in this case very utopic (see the previous example of dimension $7$), we shall determine a reduced family containing only the algebras with this characteristic sequence. For this it is sufficient to compute the $2$-cocycle of the Chevalley cohomology of $\g_{1,1,n-5}$ satisfying
\begin{equation}
\label{3step}
\left\{
\begin{array}{l}
\delta_{C,\mu_0}^2(\varphi)=0,\\
\varphi \bullet \varphi =0,\\
\delta_{R,\mu_0}^2(\varphi)=0,\\
\mu_0 \circ_1\varphi ^2+\varphi\circ_1 \mu_0 \circ_1 \varphi  +\varphi^2\circ_1 \mu_0  =0, \\
\varphi  ^3=0,
\end{array}
\right.
\end{equation}
with 
$\delta_{R,\mu_0}^2(\varphi)=\mu_0^2 \circ_1\varphi +  \mu_0  \circ_1 \varphi \circ_1 \mu_0  +\varphi\circ_1 \mu_0^2   =0.$
But any linear deformation of $\mu_0$ associated with a such cocycle belong to $\mathcal{F}^{n,3}$ that is the family of $3$-step nilpotent Lie algebras and not necessarily in the sub-family $\mathcal{F}^{n,3}_{1,1,n-5}$ of $3$-step nilpotent Lie algebras with characteristic sequence $(3,2,1,\cdots,1)$.

\begin{lemma}
We have
$$\mathcal{F}^{n,3}_{1,1,n-5}=\{\mu \in \mathcal{F}^{n,3}, \ {\rm with} \ {\rm breadth}(\mu)=3\}.$$
\end{lemma}
From the construction of the linear deformation between any Lie algebra of $\mathcal{F}^{n,3}_{1,1,n-5}$ and $\mu_0$, we can assume that $\varphi(X_1,X)=0$ for any $X$. Now, since $\mu=\mu_0+ \varphi \in \mathcal{F}^{n,3}_{1,1,n-5}$, we have necessarily $$\varphi(X_i,X_j) \in \K\{X_3,X_4,X_6\}.$$
The previous  identities and lemma imply that $\varphi$ satisfies:

\begin{equation}
\label{3stepvarphi}
\left\{
\begin{array}{l}
\varphi(X_2,X_3)=c_1X_4+e_1X_6+f_{23}X_7,\\
\varphi(X_2,X_5)=c_2X_4+e_2X_6+f_{25}X_7,\\
\varphi(X_2,X_6)=c_3X_4+e_3X_6+f_{26}X_7,\\
\varphi(X_3,X_5)=(-c_3+b_2)X_4+(d_2-e_3X_6)+f_{26 }X_7,\\
\varphi(X_5,X_6)=c_4X_4+f_{56}X_7,\\
\varphi(X_2,X_i)=c_{2i}X_4+e_{2i}X_6,  \ i \geq 7,\\
\varphi(X_5,X_i)=b_{5i}X_3+c_{5i}X_4+e_{5i}X_6,  \ i \geq 7,\\
\varphi(X_6,X_i)=b_{5i}X_4
\end{array}
\right.
\end{equation}

\end{document}